# Note on Integer Factoring Methods I

N. A. Carella, February, 2006.

***Abstract.*** This note presents the basic mathematical structure of a new integer factorization method based on systems of linear Diophantine equations. The estimated theoretical running time complexities of the corresponding algorithms are encouraging and improve the current ones. The work is presented as a theoretical contribution to the theory of integer factorization.

## 1 Introduction

This note presents the basic mathematical structure of a new integer factorization method. The new method based on systems of linear Diophantine equations facilitates the development of various integer factoring algorithms. Three integer factoring algorithms will be described, probabilistic and deterministic. These are Euclidean class algorithms and appear to be practical. The estimated theoretical running time complexities of these algorithms are favorable and improve the current running time complexities. The work is presented as a theoretical contribution to the theory of integer factorization.

The main result Theorem 1 and the corresponding algorithms are presented in section 2. The descriptions of the result and the algorithms are almost self-contained. An application to Cryptography is given in Corollary 1. The basic technical background and other relevant results (including new ones) of independent interest are given in the Appendix.

## 2 Linear Method

Let $N = pq$ be a balanced composite integer such that $p < q < 2p$ and let $\varphi(N) = N + 1 - p - q = N + 1 - t$. The sum of the prime factors of a balanced composite $N$ satisfies the inequalities $2\sqrt{N} < p + q < 1.5\sqrt{2N}$. Further, assume that $r, s = \varphi(N)$ are relatively prime integers, and consider the linear Diophantine equation

$$ru + sx = c_0, \qquad (1)$$

where $u$, $s$, and $x$ are unknown, and $c_0$ is a constant. If it has a solution $(u_0, x_0)$ such that $0 \leq |u_0|$, $|x_0| < N^{0.3}$, and half of the digits of $s$ are known, then the equation is solvable using the continued fraction algorithm or lattice reduction techniques applied to $x(N + 1 - t) \equiv c_0 \bmod r$.

The special case (the public key equation) $de - k\varphi(N) = 1$ with $0 \leq e < N^{1/4}$, $0 \leq k < N^{1/4}$, was solved by [W] using the continued fraction algorithm, see [D], [DW], [VT] for the current and improved analysis. After about fifteen years of heavy research several authors have developed (heuristic) lattice reduction techniques to improve the range of the variables $u$, $x$ in equation (1) to approximately $|u|, |x| < N^{0.3}$, see [BH], [H], [BM] etc., for the most recent developments and the references therein.

Three approaches will be considered in the development of new integer factoring algorithms based on systems of linear Diophantine equations. Toward this end, consider the system of Diophantine equations

$$ru + sx = c_0 \qquad\qquad (2)$$
$$rv + sy = c_1,$$

where $s$, $u$, $v$, $x$, and $y$ are unknown, and $c_0$, and $c_1$ are constants. Eliminating one variable, and reducing modulo $r$ lead to

$$c_0^{-1} c_1 x \equiv y \bmod r. \qquad\qquad (3)$$

Contingent on the known data and the size of the solutions, the techniques used to solve the system of equations (2), see [MK], [KS], [M], etc, are simpler than the techniques used to solve a single equation (1). This is linear techniques versus nonlinear techniques.

### Probabilistic Algorithm

An effective procedure for selecting the constants $c_0$ and $c_1$ in (2) with $|xy| < r$ is at the heart of the new integer factoring method, see Theorems 22 and 23 in the Appendix for technical details. Such procedure together with Thue's algorithm, see Theorems 12 and 13, could emerge as a very effective integer factoring algorithm. In contrast, the current best (heuristic) probabilistic algorithms have subexponential time complexities, see [LA], [CP].

***Theorem* 1.** Suppose there is an effective procedure for selecting the constants $c_0$ and $c_1$ in (2) such that $|xy| < r$ with nonnegligible probability $O(\log(N)^{-4B})$, some constant $B > 0$. Then the integer $N$ can be factored in probabilistic polynomial time.

Proof: Without loss in generality assume that the integer $N = pq$ has balanced prime factors $p < q < 2p$, and let $r$, $s = \varphi(N)$ be relatively prime integers with $r > s$. The hypothesis on the time complexity of selecting $c_0$ and $c_1$ implies that equation (3) can be solved for $x$ and $y$ in probabilistic polynomial time. Thus $\varphi(N) \equiv c_0 x^{-1} \equiv c_1 y^{-1} \bmod r$, this is an exact integer solution whenever $r > \varphi(N)$. Accordingly, the prime factors of $N$ are the roots of the equation $X^2 - (N + 1 - \varphi(N))X + N = 0$. The complexity of the algorithm is $O((\log N)^c)$ arithmetic operations on integers of size $O(\log N)$, some constant $c \leq 8$. ∎

In general, the size $|x|, |y| < r = \max\{r, s\}$ of a solution $(x, y)$ of equation (2) with $0 \leq |c_0|, |c_1| < rs < r^2$ is a random vector in $\mathbb{Z}_r \times \mathbb{Z}_r$. The solutions near $c_0^{-1}c_1 / r$ can be determined by means of either Thue's algorithm, which handles solutions of size $|xy| < r$, or the extended version called the *linear congruence algorithm*, which handles solutions of wider range, see Theorems 12, and 13 in the Appendix for more information.

In the application to the public key equation, the system of equations (2) becomes

$$eu + \varphi(N)x = 1 \qquad\qquad (4)$$
$$ev + \varphi(N)y = c_1,$$

where $r = e$ is known, and $d = u$, and $s = \varphi(N) > r$ are unknown. This requires minor changes in the calculation of $\varphi(N) = ar \pm (c_0 x^{-1} \bmod r)$, some small integer $a \geq 1$. The given data ($r = e$, $c_0 = 1$ and the estimate $|x| \leq N^\delta$, $0 < \delta < 1$) improves the chances of finding a solution since only one constant $c_1$ has to be selected at random.

As stated before, for a single equation $eu + \varphi(N)x = 1$, the current algorithms can determine the solutions up to about $|u|, |x| < N^{0.3} \leq N^\delta$, and there is some speculation that $\delta < 1/2$ is the limit, see [BH], [DW]. In contrast, the new technique using the system of equations (4) has excellent theoretical performance within the range $|x|, |y| \leq N^\delta$, $\delta \leq 1/2$, and it probably has no limit on the parameter $\delta < 1$.

***Corollary* 1**. Let $d \leq N^{1/2}$. Assume that a sequence of numbers $c_1 = c_1(t) = ev + \varphi(N)y$, where $|y| \leq N^{1/2}$ for some $0 \leq t \leq O(\log(N)^{2B})$, can be generated with nonnegligible probability $O(\log(N)^{-2B})$, some constant $B > 0$. Then the encryption key $d$ can be determined in probabilistic polynomial time.

Proof: The hypothesis $d \leq N^{1/2}$ implies that $e \in (N^{1/2}, N - N^{1/2})$. In fact $e = (k\varphi(N) + 1)/d$, where $k \equiv -\varphi(N)^{-1} \bmod d$ is a uniform random variable in the interval $[1, d-1]$, see Theorem 15. These data lead to the range of values $1 \leq de/\varphi(N) \leq N^{1/2}$, which gives $|k| < d$. See Theorem 24 for a similar result and other details. ∎

Note that the constraint $|e - \varphi(N)| = N^\theta$, $0 \leq \theta < 1$, implies that the limit $\varphi(N)/e \to 1$ as $N$ increases, for a random $d$ this condition is expected to hold. This in turn means that $|k| \approx |d|$, in fact, $|k| < |d|$, see Theorem 16. In contrast, the selection of small $e$ implies that the limit $e/\varphi(N) \to 0$, and $|k| \to 0$, $|d| \to N$ as $N$ increases.

This result seems to confirm the speculation that integer factorization is harder than determining the decryption key $d$ of the public key equation $de - k\varphi(N) = 1$, but not exponentially harder.

**Algorithm I**

Input: $N = pq$.

Output: $p, q$.

1. Put $m = 0$, and select a random integer $r \in [N - 1.5\sqrt{2N}, N]$.

2. While $m \leq 4(\log N)^4$, select a pair of random integers $c_0, c_1 \in [1, rs]$, $rs \approx (N - 1.5\sqrt{2N})^2$.

3. Compute the solutions $(x_i, y_i)$ of $c_0^{-1}c_1 x \equiv y \bmod r$ near $c_0^{-1}c_1 / r$, for $0 \leq i < 5\log r$.

4. Compute $s_{0,j} \equiv c_0 x_j^{-1} \bmod r$, $s_{1,j} \equiv c_1 x_j^{-1} \bmod r$, $s_{2,j} \equiv c_0 y_j^{-1} \bmod r$, $s_{3,j} \equiv c_1 y_j^{-1} \bmod r$, and $D_{i,j} = (N + 1 \pm s_{i,j})^2 - 4N$, for $0 \leq i < 4$, $0 \leq j < 5\log r$.

5. If $D_{i,j}$ is a square, then compute the roots $X_{i,j}$ $Y_{i,j}$ of $X^2 - (N + 1 \pm s_{i,j})X + N = 0$.

6. If $X_{i,j}Y_{i,j} \neq \pm N$ then repeat Step 2.

7. Return $p = |X_{i,j}|$, $q = N/p$ or $p = |Y_{i,j}|$, $q = N/p$.

In step 1, the condition $\gcd(r, s) = 1$ has non negligible probability, see Theorem 25. In step 2, a random selection of the constants is one of many possible strategies. In step 5, a search for a square among the discriminants $D_{i,j} = (N + 1 \pm s_{i,j})^2 - 4N$ is conducted, and then compute the corresponding roots. The complementary discriminants $E_{i,j} = (N + 1 \pm r \pm s_{i,j})^2 - 4N$ should also be included to expedite the factorization. A simple example is provided to demonstrate the basic procedure of Algorithm I. It also suggests that the system (4) has deeper structure and the algorithm could be more flexible than predicted, see Theorem 17 for some explanation.

***Example* 1.** Given the data $N = 2257$, $e = 2431$, and the estimate $|k| < N^{3/4}$, find the decryption key $d$ in $de - k\varphi(N) = 1$.

The supplied data is used to factor $N$, and then compute $d$. Put $r = e$, $x = k$, $s = \varphi(N)$, and $c_0 = 1$ as in (4). Algorithm I (step 3 is based on Thue's algorithm) was used to probe the sequence of integers

$$c_1 = (N + 2)^2 - (N + 1)^2, \ldots, c_1 = (N + 5)^2 - (N + 1)^2 = 18080, \ldots.$$

The result for the fourth one $c_1 = 18080$, and the corresponding linear equation $1063x \equiv y \bmod r$ is shown in the table. The value $271 \equiv s_{i,j} \equiv c_i x_j^{-1} \bmod r$ gives a square discriminant $E_{i,j} = (N + 1 - r + s_{i,j})^2 - 4N = 576$, so $X^2 - (N + 1 - r + s_{i,j})X + N = 0$ has integral roots $X = 37, 61$, which are the factors of $N$. The key is $d = (305 \cdot 36 \cdot 60 + 1)/2431 = 271$.

| $I$ | $x_i$ | $y_i$ | $c_0 x^{-1} \bmod r$ | $c_1 x^{-1} \bmod r$ | $c_0 y^{-1} \bmod r$ | $c_1 y^{-1} \bmod r$ |
|---|---|---|---|---|---|---|
| 1 | 1 | 1063 | 1 | 1063 | 1889 | 1 |
| 2 | –2 | 305 | 1215 | 684 | 271 | 1215 |
| 3 | 7 | 2 | 1042 | 1541 | 1659 | 1042 |
| 4 | –16 | 16 | 2279 | 1301 | 2161 | 2279 |
| 5 | 263 | 2 | 2089 | 1104 | 608 | 2089 |
| 6 | –542 | 4 | 1063 | 1985 | 1 | 1063 |

Assuming the correctness of Theorem 23, certain sequences $\{ c_1 = c_1(t) : 0 \le t < 2(\log N)^2 \}$ of length $2(\log 2257)^2 \approx 23$ produces a solution of equation (4) and the factorization of $N = 2257$ with nonnegligible probability. This example only requires 4 trials from the set of 23 trials.

### *Approximation Algorithm*

For each fixed pair $c_0$, $c_1$, the algorithm utilized to solve the linear equation $c_0^{-1} c_1 x \equiv y \bmod r$ generates a series of approximations $(x_0, y_0)$, $(x_1, y_1)$, ..., $(x_k, y_k)$ of a solution $(x, y)$ closed to $c_0^{-1} c_1 / r$. These in turn give a series of approximations $s_{0,1}, s_{0,2}, \ldots, s_{3,k}$ of $\varphi(N)$. It is very likely that some of these approximations have sufficiently small errors $| s_{i,j} - \varphi(N) |$.

An integer $s_{i,j} \in \mathbb{Z}$ in the range $\varphi(N) - N^{1/4}$ to $\varphi(N) + N^{1/4}$ will be called an *effective approximation* of $\varphi(N)$. The effective approximations of $\varphi(N)$, where $N = pq$ is a balanced composite integer with $p < q < 2p$, are detected by means of the inequalities

$$| N - X_0 X_1 | < 3N^{3/4} \quad \text{and} \quad | X_0 | < | X_1 | < 2 | X_0 | \tag{5}$$

or similar inequalities. Here the real numbers $X_0$, $X_1$ are the roots of $X^2 - (N + 1 \pm s_{i,j}) X + N = 0$ or its complementary equation $X^2 - (N + 1 - r \pm s_{i,j}) X + N = 0$.

Since effective approximations $s_{i,j} = \varphi(N) + t$, where $0 < | t | \le N^{1/4}$, are far more frequent than exact solutions $s_{i,j} = \varphi(N)$, it prompts the idea that the integer factorization problem should be approached from an approximation point of view, and Algorithm I should be upgraded to support this improvement. Toward this goal, let $p_0 = [X_0]$, and $q_0 = [X_1]$, and suppose that inequalities (5) are satisfied. Then

$$(p_0 + U)(q_0 + V) = N \tag{6}$$

has an integral solution $(U_0, V_0)$ with $| U_0 |$, $| V_0 | < N^{1/4}$ or $| U_0 V_0 | < N^{1/2}$. In the case $U_0 = 0$ or $V_0 = 0$ there is an exact solution. And in the case $U_0 \ne 0$ or $V_0 \ne 0$, the small integral solution is resolved in deterministic polynomial time using lattice reduction techniques, see [CR].

### Algorithm II
Input: $N = pq$.
Output: $p$, $q$.
1. Put $m = 0$, and select a random integer $r \in [N - 1.5\sqrt{2N}, N]$.
2. While $m \le 4(\log N)^4$, select a pair of random integers $c_0$, $c_1 \in [1, rs]$, $rs \approx (N - 1.5\sqrt{2N})^2$.
3. Compute the solutions $(x_i, y_i)$ of $c_0^{-1} c_1 x \equiv y \bmod r$ near $c_0^{-1} c_1 / r$, for $0 \le i < 5 \log r$.
4. Compute $s_{0,j} \equiv c_0 x_j^{-1} \bmod r$, $s_{1,j} \equiv c_1 x_j^{-1} \bmod r$, $s_{2,j} \equiv c_0 y_j^{-1} \bmod r$, $s_{3,j} \equiv c_1 y_j^{-1} \bmod r$, and $t_{i,j} = r - s_{i,j}$, for $0 \le j < 5 \log r$.
5. Find an effective approximation among the integers $s_{i,j}$, and $t_{i,j} = r \pm s_{i,j}$ using the roots $X_{i,j}$ $Y_{i,j}$ of $X^2 - (N + 1 \pm s_{i,j}) X + N = 0$, or $X^2 - (N + 1 \pm t_{i,j}) X + N = 0$, for $0 \le i < 4$, $0 \le j < 5 \log r$.

6. Put $p_0 = [X_{i,j}]$, $q_0 = [Y_{i,j}]$ and compute the integral roots of $UV + q_0U + p_0V + p_0q_0 - N = 0$.

7. Return $p = |\,p_0 + U_0\,|$, $q = N/p$.

### *Deterministic Algorithm*

The best deterministic rigorously analyzed integer factorization algorithm has exponential time complexity of $O(N^{1/4})$, see [CP], [MP], and earlier works in this direction [LN], [LA], et cetera. This section supplies the preliminary details of the new technique that have the potential to improve the current time complexity.

Let $p_n/q_n$ be a convergent of $r/s$, and for $n \geq 2$, write

$$rq_n - sp_n = z_n$$
$$rq_{n+1} - sp_{n+1} = z_{n+1}. \tag{7}$$

Canceling the variable $s$ leads to $(p_nq_{n+1} - p_{n+1}q_n)r + z_np_{n+1} - z_{n+1}p_n = 0$. Another round of simplification using the properties of convergents (Theorem 6) yields

$$z_np_{n+1} - z_{n+1}p_n = (-1)^n r. \tag{8}$$

### Algorithm III

Input: $N = pq$.

Output: $p$, $q$.

1. Select a random pair $r, s \in [N - 2\sqrt{N}, N - 1.5\sqrt{2N}]$, $\gcd(r, s) = 1$, and compute the convergents $p_n/q_n$ of $r/s$.

2. Compute a solution $(z_n, z_{n+1})$ of $z_np_{n+1} - z_{n+1}p_n = (-1)^n r$, for some $1 \leq n < 5\log r$.

3. Compute $t = (z_n + (N+1)p_n - rq_n)/p_n$, and the roots $X_0, X_1$ of $X^2 - tX + N = 0$.

4. If $x_0x_1 \neq \pm N$, then repeat Step 2.

5. Return $p = |\,X_0\,|$, $q = N/p$.

In step 1, the balanced property provides a concrete range to select a pair of integers $r, s \in [N - 2\sqrt{N}, N - 1.5\sqrt{2N}]$ with $\gcd(r, s) = 1$, where $s$ is an estimate of $\varphi(N)$. In step 2, the solution $z_n, z_{n+1}$ can be obtained by means of linear congruence, Euclidean algorithm or other techniques. This is a deterministic process. The main obstacle is the determination of nontrivial solution $z_n, z_{n+1}$ of equation (8).

Here the inequalities

$$q_n < N^\alpha \leq q_{n+1} \quad \text{and} \quad |\,rq_n - sp_n\,| = z_n \leq N/q_{n+1} \leq N^{1-\alpha}, \tag{9}$$

where $0 < \alpha < 1$, can serve as a guide in the determination/minimization of $z_n$, $z_{n+1}$. These inequalities are derived from the well known convergent's inequality

$$\left| \frac{r}{s} - \frac{p_n}{q_n} \right| \leq \frac{1}{q_n q_{n+1}} . \tag{10}$$

The minimization of $z_n$, $z_{n+1}$ plays against the maximal value of the convergents $p_n/q_n$ computable by means of the approximation $r/s$ of $r/\varphi(N)$. Some works on the maximal computable convergents $p_n/q_n$ in polynomial time appear in [VT], and [D].

**Research Problems**

1. Prove that the construction of an effective approximation $s_{i,j}$ of $\varphi(N)$ has deterministic polynomial time complexity.
2. Determine the distribution of pairs $c_0$, $c_1$ that give effective approximations $s_{i,j}$ of $\varphi(N)$.

## 3 Appendix

A limited introduction to the theoretical foundation of the new integer factorization method is supplied here, extensive technical details on these subjects are available in the literature. A few auxiliary concepts of interest in future applications and developments are also included.

**Continued Fractions**

The greatest common divisor $d = \gcd(a, b)$ of a pair of integers $a, b \in \mathbb{N}$ is the greatest integer $d \in \mathbb{N}$ that divides both.

*Euclidean algorithm*: Given a pair of integers $a, b > 0$, the algorithm computes the greatest common divisor $\gcd(a, b) = r_{n+k-1}$ by the simple process

$$r_0 = a, \qquad r_1 = b, \qquad r_n = a_{n+1}r_{n+1} + r_{n+2}, \qquad 0 < r_{n+1} < r_n,$$

where $n \geq 2$, and $r_{n+k} = 0$ for some $k \geq 0$.

The decreasing sequence $r_0 > r_1 > \cdots > 0$ ensures that it is a finite process. This follows from the recursive formula $r_n = a_{n+1}r_{n+1} + r_{n+2} \geq r_{n+1} + r_{n+2} > 0$, and the well ordered principle: every nonempty set of nonnegative integers contains a smallest integer.

The Euclidean algorithm is one of the oldest mathematical results, and still is a topic of current research. This algorithm is documented in Euclid's Elements, circa 300 BC, but was probably discovered earlier. This basic result is a sine qua non of computational number theory.

***Theorem* 2.** (Lame 1844) For $a, b \leq N$, the Euclidean algorithm has the complexity of $O(\log N)$ integer divisions.

**Theorem 3.** (HL) The number of integer divisions performed by the Euclidean algorithm on the set of pairs $1 \leq a$, $b \leq N$ is normally distribute with mean $\mu \approx ((12\log 2)/\pi^2)\log(N)$ and variance $\sigma \approx c\log(N)$ as $N \to \infty$, where $c > 0$ is a constant.

The continued fraction of a real number $\xi \in \mathbb{R}$ is defined by the sequence $\xi = [a_0, a_1, a_2, \ldots] \in \mathbb{Z} \times \mathbb{N}^\infty$. The inverse is simply $\xi^{-1} = [0, a_0, a_1, a_2, \ldots]$. Each rational number has a unique finite continued fraction with $a_i = [r_i / r_{i+1}]$ and the last quotients $a_n = 1$ or $a_n = 2$. Otherwise for irrational number the continued fraction is infinite with

$$a_0 = [\xi_0], \xi_0 = \xi, \qquad a_1 = [\xi_1]^{-1}, \xi_1 = \xi_0 - a_0, \qquad a_2 = [\xi_2]^{-1}, \xi_2 = \xi_1 - a_1, \ldots .$$

The continued fractions of real numbers are representations of the numbers which have different properties other than those properties of decimal expansions of real numbers. One strikingly distinct property is the distribution of the digits. The distribution of the $n$th digit $d_n$ of the decimal expansions of real numbers have uniform distribution with probability density function $P(d_n = k) = 1/10$, but the distribution of the digits of the continued fraction expansions of real numbers do not have uniform distribution. It was conjectured by Gauss that the probability density function $P(a_n = k)$ of the $n$th digit of the continued fraction of real numbers converges to $P(a_n = k) = \log_2(1 + 1/k(k+2))$ as $n$ tends to infinity.

**Theorem 4.** (Gauss-Kuzmin 1928) Let the sequence $[a_0, a_1, a_2, \ldots] \in \mathbb{Z} \times \mathbb{N}^\infty$ be the continued fraction of a real number. Then the probability that the $n$th digit $a_n = k \geq 1$ is $\log_2(1 + 1/k(k+2))$ as $n \to \infty$.

The probability that the first digit $a_1 = k$ in the continued fraction of $\xi \in (0, 1)$ is given by $P(a_1 = k) = 1/k(k+1)$ because $\xi \in (1/(k+1), 1/k)$. For example,

$$P(a_1 = 1) = 1/2, \qquad P(a_1 = 2) = 1/6, \qquad P(a_1 = 3) = 1/12, \ldots .$$

The probability $P(a_1 = k)$ of the first digit $a_1$ is independent, but the probabilities $P(a_n = k)$ of the $n$th digits $a_n$ with $n \geq 2$ are somewhat dependent and nearly independent for large $n$, consult the literature for advanced details.

**Theorem 5.** ([HW]) The $n$th digit $a_n$ is unbounded for almost all real numbers, the set of real numbers which have bounded $a_n \leq k$, $k$ constant, is a null set.

The $n$th convergent of the continued fraction is given by the ratio $p_n/q_n$, where

$$p_0 = 0, \qquad p_1 = 1, \qquad p_n = a_{n-1}p_{n-1} + p_{n-2},$$
$$q_0 = 1, \qquad q_1 = a_1, \qquad q_n = a_{n-1}q_{n-1} + q_{n-2}, \qquad n \geq 2.$$

**Theorem 6.** The convergents satisfy the followings relations.
(i) $p_nq_{n-1} - p_{n-1}q_n = (-1)^{n-1}$,        (ii) $p_nq_{n-2} - p_{n-2}q_n = (-1)^n a_n$,
(iii) $\gcd(p_n, q_n) = 1$, and $\gcd(p_n, p_{n+1}) = 1$.

***Theorem 7.*** (Dirichlet 1842)  Let $p, q \geq 0$ be integers such that $\gcd(p, q) = 1$. If the inequality

$$\left| \xi - \frac{p}{q} \right| < \frac{1}{2q^2} \tag{11}$$

holds, then $p/q = p_n/q_n$ is a convergent of the continued fraction of $\xi$.

## Linear Diophantine Equations

***Theorem 8.*** Let $r_1, r_2, \ldots, r_k \in \mathbb{Z}$ be integers such that $\gcd(r_1, r_2, \ldots, r_k) = d$. Then $r_1 x_1 + \cdots + r_k x_k = n$ has an integer solution if and only if $d$ divides $n$.

***Theorem 9.*** Let $r, s \in \mathbb{N}$ be relatively prime integers. Then $rx + sy = 1$ is solvable.

Proof: Compute the penultimate convergent of $r/s$, and put $x = q_{n-1}$ and $y = -p_{n-1}$.  ∎

It is clear that there is a unique solution of $rx + sy = 1$ such that $0 < x < s$. Furthermore, the ratio $r/s$ gives a very precise estimate $y \approx -(r/s)x$.

***Theorem 10.*** (Ariabhata 499)  Let $r, s \in \mathbb{N}$ be relatively prime integers. Then the solutions of $rx + sy = n$ are given by $x = x_0 + st$, $y = y_0 - rt$, $t \in \mathbb{Z}$.

The initial solution $(x_0, y_0)$ is constructed using the previous result or similar techniques.

***Theorem 11.*** (Bachet 1612?)  The linear congruence $ax \equiv b \bmod N$ is solvable if and only if $\gcd(a, N) = d$ divides $b$. A solvable equation has $d$ solutions $x = x_0 + (N/d)t$, $0 \leq t < d$.

Proof: It follows from Theorem 10.  ∎

***Theorem 12.*** (Thue 1900)  Let $N > a$, and $\gcd(a, N) = 1$. Then $x = r_k$, $y = q_k$ is a solution of $ax \equiv y \bmod N$, with $|x|, |y| \leq \sqrt{N}$.

Proof: Write $ax = y + zN$, and assume that $0 < xy < N$. Then

$$\left| \frac{a}{N} - \frac{z}{x} \right| = \left| \frac{y}{xN} \right| < \frac{1}{x^2} . \tag{12}$$

Here $z/x = p_n/q_n$ is the $n$th convergent of the continued fraction of $a/N$, and $r_0 = a$, $r_1 = N$, $r_n = a_{n+1}r_{n+1} + r_{n+2}$, $q_n = a_{n-1}q_{n-1} + q_{n-2}$, $n \geq 2$. The rest follows from Dirichlet's Theorem.  ∎

A related proof and applications appear in [M]. A recent work [KS] improves the range of $x$ and $y$ in $ax \equiv y \bmod N$ to an arbitrary size $B < N$, but it seems to depend on $a$.

**Theorem 13.** ([KS]) Given integers $a, c, N, B$ where $0 < a, B < N$, $0 \leq c < N$, and $\gcd(a, N) = 1$ holds, the *Linear Congruence Algorithm* determines the unique minimal solution $0 \leq x < N$, $0 \leq y < B$ of the congruence $ax \equiv y + c \bmod N$.

To apply this algorithm to the equation $ax \equiv y \bmod N$ with $c = 0$, a linear change of variables is required, the authors recommend $(x, y) \rightarrow (x - 1, y)$. In this case, nontrivial solutions $(x, y) \neq (0, 0)$ requires the parameter $B \leq a$. The basic procedure of the linear congruence algorithm is given below, the reader should consult [KS] for a complete analysis.

*Linear Congruence Algorithm*
Input: $a, c, N, B$, where $0 < a, B < N$, $0 \leq c < N$, and $\gcd(a, N) = 1$.
Output: $x_0, y_0$ such that $ax_0 = y_0 + c \bmod N$ and $x_0 \geq 0$ is minimal such that $0 \leq y_0 < B$.
1. Set $a' = a$, $c' = c$, $N' = N$
2. Set $y' = -c' \bmod N'$
3. While $y' \geq B$ do
4. Set $(a', N') = (-N' \bmod a', a')$
5. Set $c' = c' \bmod N'$, $y' = -c' \bmod N'$
6. Set $y_0 = y'$, $x_0 = a^{-1}(y_0 + c) \bmod N$
7. Return $(x_0, y_0)$.

**The Inverse of a Modulo $N$**
Let $x \in [1, N)$ be relatively prime to $N$, and let $x^{-1}$ be the inverse of $x$ modulo $N$. It is known that inverses are uniformly distributed in the interval $[1, N)$. The formal description of this phenomenon is as follows. Consider the map

$$f_N(x) = \left( \left\langle \frac{x}{N} \right\rangle, \left\langle \frac{x^{-1}}{N} \right\rangle \right), \tag{13}$$

where $< z > \in [0, 1]$ is the fractional part of the real number $z \in \mathbb{R}$. The function $f_N$ maps a random point $x \in [1, N)$ such that $\gcd(x, N) = 1$ to a random pair $f_N(x) = (a, b) \in [0, 1]^2$. The pairs are uniformly distributed in the unit rectangle $[0, 1]^2$.

**Theorem 14.** ([BK]) Let $R \subset [0, 1]^2$ be a measurable set having the property that for every $\varepsilon > 0$ there exists a finite collection of nonoverlapping rectangles $R_1, R_2, \ldots, R_k$, such that $\cup R_i \subseteq R$ and the area $R - \cup R_i < \varepsilon$. Then

$$\lim_{N \to \infty} \frac{\# \mathrm{Im}(f_N) \cap R}{\varphi(N)} = \mathrm{Area}(R). \tag{14}$$

Here the term $\mathrm{Im}(f_N) = \{ f_N(x) = (a, b) \in [0, 1]^2 : x \in [1, N) \text{ and } \gcd(x, N) = 1 \}$. The proof of this interesting result is based mostly on exponential sums analysis.

An elementary analysis will be utilized to obtain information on the distribution of the inverses modulo $N$. This analysis starts with the simple observation that the pattern of the inverses of small integers is the following:

(i) $2^{-1} \equiv (\vartheta N - 1)/2 + 1 \bmod N$, with $N \not\equiv 0 \bmod 2$, for some $0 < \vartheta < 2$,

(ii) $3^{-1} \equiv (\vartheta N - 2)/3 + 1 \bmod N$, with $N \not\equiv 0 \bmod 3$, for some $0 < \vartheta < 3$,

(iii) $4^{-1} \equiv (\vartheta N - 3)/4 + 1 \bmod N$, with $N \not\equiv 0 \bmod 2$, for some $0 < \vartheta < 4$,

and similar expressions for other integers. This simple observation provides a much more precise characterization on the form and the distribution of the inverses.

***Theorem* 15.** Let $a$, $N$ be a pair of integers. Then the followings hold.
(i) If $a$ is relative prime to $N$, then the inverse $a^{-1} \bmod N$ is precisely the integer

$$a^{-1} = \frac{(-N^{-1} \bmod a)N + 1}{a} \in \mathbb{N}. \tag{15}$$

(ii) The integer $\vartheta \equiv -N^{-1} \bmod a$ is a uniform random variable in the interval $[1, a - 1]$.

The proof is simple, it just requires a bit of tinkering with the calculation of inverse. That is, rewrite $ax \equiv 1 \bmod N$ as

$$x \equiv \frac{\vartheta N - (a - 1)}{a} + 1 \bmod N, \quad \text{some } 0 < \vartheta < a. \tag{16}$$

Now use elementary manipulations and the fact that the inverse is unique to assemble it. By symmetry the same result should also holds for $x$ in place of $a$. The uniform property of the integer $\vartheta \in [1, a)$ follows from Theorem 14.

As a spin-off, the new formula (15) (or its counterpart in other finite rings) appears to be the fastest method around for computing the inverses of small $| a | \leq O(\log(N)^c)$ residues modulo $N$, it is faster than the extended Euclidean algorithm $ax + Ny = 1$ or the exponentiation method $a^{-1} \equiv a^{\varphi(N)-1} \bmod N$, (Fermat's little Theorem). The inverse can also be computed using the (Voroni 1900) formula

$$a^{-1} \equiv (3 - 2a + 6 \sum_{k=1}^{a-1} [kN/a]^2) \bmod N, \tag{17}$$

but this is not as effective as formula (15) for any $a < N$. The complexity is $O(\log\log N)$ divisions versus $O(\log N)$.

For fixed $N$, the size $| a |$ of a residue $a \in \mathbb{Z}_N$ serves as a measure of the randomness of the inverse $a^{-1}$. The randomness increases as $| a |$ increases to $N/2$, but it decreases as $| a |$ decreases to 1, see statement (ii) above. Accordingly, the inverse of an integer of size $a = N^\alpha < N$ must be in the range

$$N^{1-\alpha} \leq x \leq N - N^{1-\alpha}, \tag{18}$$

where $0 < \alpha < 1$. In particular,

If $a \in (1, N^{1/4})$, then the inverse is $a^{-1} \in (N^{3/4}, N - N^{3/4})$,
If $a \in (1, N^{1/2})$, then the inverse is $a^{-1} \in (N^{1/2}, N - N^{1/2})$,
If $a \in (N^{1/4}, N^{3/4})$, then the inverse is $a^{-1} \in (N^{1/4}, N - N^{1/4})$,

and so on. This implies that as $a$ varies over the interval $(1, N^{1/2})$ the inverse $a^{-1}$ is not uniformly distributed in $(1, N)$. Moreover, it seems that the inverses of a small number of points in $\mathbb{Z}_N$ have nonrandom inverses. For example, the inverse of 5 modulo $N$ is known a priori for any $N \not\equiv 0 \bmod 5$, id est,

$$5^{-1} \equiv \begin{cases} (4N-4)/5+1 & \text{if } N \equiv 1 \bmod 5, \\ (2N-4)/5+1 & \text{if } N \equiv 2 \bmod 5, \\ (3N-4)/5+1 & \text{if } N \equiv 3 \bmod 5, \\ (N-4)/5+1 & \text{if } N \equiv 4 \bmod 5. \end{cases} \tag{19}$$

The following results establish relationships between the distance between a pair of integers and the solutions of the corresponding linear Diophantine equation.

**Theorem 16.** Let $rx + sy = n$, where $r, s \in \mathbb{N}$ are fixed relatively prime integers. Then
(i) $x \equiv nr^{-1} \bmod y$ is a uniform random variable in $[1, y)$.
(ii) $y \equiv ns^{-1} \bmod x$ is a uniform random variable in $[1, x)$.

Proof: It follows from either Theorem 14 or 15. ∎

For a fixed pair $r, s \in \mathbb{N}$, the magnitude $|x|$ of the integer $x \in \mathbb{Z}$ with respect to $y$ depends on the distance $|r - s|$. In particular, if $|r - s| = o(r)$, then $|x| < |y|$.

**Theorem 17.** Let $r, s \in \mathbb{N}$ be relatively prime integers, and let $r - s = t$. Then the solutions of $rx + sy = n$ satisfy the congruences $xt \equiv n \bmod s$ and $yt \equiv -n \bmod r$.

Proof: Consider the system of equations $rx + sy = n$ and $r - s = t$. ∎

These results have interesting applications to the public key equation $de - k\varphi(N) = 1$. Specifically, if $e \approx N$, and the parameter $k = N^\delta < N$, $\delta > 0$, then the distance $|e - \varphi(N)| \in (e^{1-\delta}, e - e^{1-\delta})$. Similarly, for $|e - \varphi(N)| = o(N)$, and the parameter $k = N^\delta$, then the decryption key is

of size $d \approx N^\delta$. In fact $|k| < |d|$, and the parameter $k \equiv -\varphi(N)^{-1}$ mod $d$ is a uniform random variable in the interval $[1, d)$.

The size of the solutions and the maximal distance of the pairs $x, y \in \mathbb{Z}$ such that $xy \equiv c$ mod $N$ are of considerable interest in the theory of linear congruences and their applications.

**Theorem 18.** ([GV]) For any fixed $\varepsilon > 0$ and any prime $p$ the set of residues $\{ xy$ mod $p : 1 \le x,$ $y \le p^{1/2}(\log p)^{2+\varepsilon} \}$ contains $(1 + o(1))p$ residues classes modulo $p$.

The more general form of this result claims that the congruence equation $xy \equiv c$ mod $N$ has a solution $1 \le x, y \le N^{1/2}(\log N)^{2+\varepsilon}$ for almost any $c < N$. This could be of considerable interest in future development of the integer factoring technique presented here.

**Theorem 19.** ([KM]) For a given $n \ge 2$ let $M(n) = \max \{ |a - b| : a, b \in \{ 1, 2, \ldots, n \}$ and $ab \equiv 1$ mod $n \}$. Then $M(n) \le [n - 2\sqrt{n-1}]$ with equality if and only if $n = m^2 + lm + 1$, and $0 \le l < 2\sqrt{m} + 1$.

The proof of this result springs from the arithmetic-geometric mean inequality $(a + b)/2 \ge \sqrt{ab}$ for nonnegative numbers $a, b \in \mathbb{R}$.

**Frobenius Problem**
The *Frobenius number* $f(r_1, r_2, \ldots, r_k)$ is the largest integer $n > 0$ that is not representable by the linear equation $n = r_1 x_1 + \cdots + r_k x_k$, where $r_1, r_2, \ldots, r_k$ are constants and $x_i \ge 0$.

In the 2 dimensional case, the representation of $rx + sy = n$ is obtained from the general solution $x = x_0 + st, y = y_0 - rt, t \in \mathbb{Z}$. If it does have a solution, it is determined by locating a value of the parameter $t$ such that $x_0 + st \ge 0$ and $y_0 - rt \ge 0$.

**Theorem 20.** Let $r, s, t \in \mathbb{N}$ be fixed integers such that $\gcd(r, s, t) = 1$. Then
(i) $f(r, s) = rs - r - s$.
(ii) $f(r, s, t) \ge \sqrt{3rst} - r - s - t$.

The origin of the first fact is unknown, [R, p. 31], and the second fact is due to [DN].

**Theorem 21.** Let $r, s \in \mathbb{N}$ be relatively prime integers, and let $R(n)$ denotes the number of representations of $n = rx + sy, x, y \ge 0$. Then it is given by the (Popoviciu 1953) formula

$$R(n) = \frac{n}{rs} - \left( \left( \frac{s^{-1}n}{r} \right) \right) - \left( \left( \frac{r^{-1}n}{s} \right) \right) + 1, \qquad (20)$$

where $((x))$ is the fractional part of $x$, and $s^{-1}$ mod $r$ is an integer in the interval $[1, r)$.

The number $R_{i,j}$ of solutions $x_1 \geq 0, \ldots, x_k \geq 0$ of a system of two equations

$$\begin{aligned} r_1 x_1 + \cdots + r_k x_k &= n, \\ s_1 x_1 + \cdots + s_k x_k &= m, \end{aligned} \qquad (21)$$

is given by the coefficient $R_{n,m}$ of $X^n Y^m$ in the expansion

$$\prod_{i=1}^{k} \left(1 - X^{r_i} Y^{s_i}\right)^{-1} = \sum_{n,m=1}^{\infty} R_{n,m} X^n Y^m. \qquad (22)$$

Counting the number of solutions of a system of equations seems to be more difficult than a single equation, however, this is an old result of unknown author, see [TL, p. 188].

**Statistical Properties of the Solutions**

***Theorem 22.*** Let $r, s \in \mathbb{N}$ be relatively prime integers such that $r > s$. Then half of the integers between $s$ and $(r-1)(s-1)$ are representable as $n = rx + sy$, some $x, y \geq 0$.

There are many proofs of this result due to Sharp 1883, see [R, p. 104] for extensive details. One interesting proof amounts to showing that either $n \geq 1$ or $rs - n$ has a representation as $rx + sy$, some $x, y \geq 0$. This is equivalent to showing that $R(n) + R(rs - n) = 1$, which is a routine application of equation (20).

A randomly selected integer $n$ from the interval $[s, (r-1)(s-1)]$ has a representation as $n = rx + sy$ such that $x, y \geq 0$, with probability $1/2$, this follows from Theorem 22. Moreover, the density of the integers that have such representations increases as the integers $n$ get closer to $(r-1)(s-1)$. Eventually, for every integer $n > (r-1)(s-1)$, the probability that it has such representation is 1, see Theorem 20.

**Small Solutions**

The distance $|r - s| \geq 1$ between the integers $r$ and $s$ is an important factor in the existence of integers $n = rx + sy < r$ with small $|x|, |y|$. For example,

if the distance is 1, then every integer $1 \leq |n|$ is of the form $n = rx + sy$ with $0 \leq |x|, |y| \leq |n|$,
if the distance is 2, then every even integer $2 \leq |n|$ is of the form $n = rx + sy$ with $0 \leq |x|, |y| \leq |n/2|$,
if the distance is $d$, then every even integer $d \leq |n|$ is of the form $n = rx + sy$ with $0 \leq |x|, |y| \leq |n/d|$,
etc. So some of the integers close to the distance $|r - s|$ should have solutions of small size $|x|$, $|y|$. This idea is the basis of step 1 in Algorithm I.

The distribution of small solutions and an effective procedure for identifying integers $n$ that have representations as $n = rx + sy$ with $0 \leq |x|, |y| \leq r^{1/2}$, or $0 \leq |xy| \leq r$, or some other specified range are important problems in the theory of linear Diophantine equations and their

applications. Below there is a sketch of an attempt to show that there is a procedure for selecting such integers in probabilistic polynomial time with non negligible probability $O(\log(N)^{-2B})$, some constant $B > 0$.

**Theorem 23.**   Let $r$ and $s$ be relatively prime integers such that $s < r - r^{1/2}\log(r)^{-B}$, $B > 0$ constant. Then the integer $n = 2r(b_2 - b_1) + b_2^2 - b_1^2$, where $0 \leq |b_i| \leq r^{1/2}\log(r)^{-B}$, has a representation as $rx + sy$, some $0 \leq x, y \leq r^{1/2}$, with probability approximately $1/(2\log(r)^{2B})$.

Proof: Let $b_2 > b_1$ be small integers, $0 \leq |b_i| \leq r^{1/2}\log(r)^{-B}$, $B > 0$. Since $r > s$, $\gcd(r, s) = 1$, and $r^2 > rs$, there exists a pair $x_i, y_i \geq 0$ such that $n_i = (r + b_i)^2 = rx_i + sy_i$ with probability 1, see Theorem 20. On the other hand, since $s < n = n_i - n_j < rs$, the integer representation

$$n = 2r(b_2 - b_1) + b_2^2 - b_1^2 = r(x_2 - x_1) + s(y_2 - y_1) = rx + sy, \tag{23}$$

where $x_i, y_i \geq 0$, exists with probability $1/2$, see Theorem 22. Moreover, if it does exist, then $0 \leq x_2 - x_1 = x \leq 6r^{1/2}\log(r)^{-B}$, and $0 \leq y_2 - y_1 = y \leq 6r^{1/2}\log(r)^{-B}$. Let

$$U = \{ \ n = 2r(b_2 - b_1) + b_2^2 - b_1^2 = rx + sy \ : 0 \leq b_i \leq r^{1/2}\log(r)^{-B} \ \},$$

and

$$V = \{ \ n = rx + sy : 0 \leq x, y \leq r^{1/2} \ \}.$$

If an integer $n$ does have such representation, then $n \in U \cap V$. Now observe that the cardinalities of these subsets are $r\log(r)^{-2B} - r^{1/2}\log(r)^{-B}$ and $r = \#V$ respectively. Therefore, the probability of finding such an integer is given by

$$\frac{\#\{ \ n \in U \cap V \ \}}{\#V} = \frac{1}{2}\frac{\#U}{\#V} \approx \frac{1}{2(\log r)^{2B}} \tag{24}$$

as $r \to \infty$.   ∎

In light of this result, it seems that the sums and differences $n = \pm b_4^2 \pm b_3^2 \pm b_2^2 \pm b_1^2$ of small squares are excellent sources of linear combinations $n = rx + sy$ with small $|x|, |y| \leq r^{1/2}$.

**Theorem 24.**   Let $r, s \in \mathbb{N}$ be relatively prime integers such that $|r - s| = r^{\theta}$, $0 \leq \theta < 1$, where $r > s$ is known and $s$ is unknown. Suppose that $ru + sx = 1$ has a solution of size $|x| < r^{1/2}$. Then $s$ can be determined in probabilistic polynomial time.

Proof: Let $B > 0$ be a constant, and construct a sequence of numbers $c = c(t) = rv + sy$ that has a small solution $|y| < r^{1/2}$ for some $0 \leq t \leq O(\log(N)^{2B})$. Now consider the system of equations

$$ru + sx = 1 \qquad (25)$$
$$rv + sy = c.$$

The hypothesis $|\, r - s \,| = r^{\theta}$, and the estimate $|\, x \,| \leq r^{1/2}$ imply that $|\, u \,| \approx |\, x \,| \leq r^{1/2}$. Thus assuming Theorems 12 and 23, it follows that a solution of the congruence $c^{-1}x \equiv y \bmod r$ can be determined in probabilistic polynomial time with probability approximately $1/(2\log(r)^{2B})$. $\blacksquare$

**Relatively Prime Integers**
Questions on the density of relatively prime integers arise on the analysis of Algorithms I, II and III. Several results for the probabilities of obtaining pairs of relatively prime integers in infinite and short intervals are considered now.

Let $N \in \mathbb{N}$. Euler function is defined by $\varphi(N) = N \prod_{\text{prime } p \mid N} (1 - 1/p)$. It enumerates the number of relatively prime integers $a < N$, and satisfies $\dfrac{N}{2} < \varphi(N) < \dfrac{N}{e^{\gamma} \log \log N}$. The universal exponent function is defined by $\lambda(N) = lcm(\lambda(p_1^{v_1}), \lambda(p_2^{v_2}), \dots, \lambda(p_k^{v_k}))$, where

$$\lambda(p^v) = \begin{cases} \varphi(2^v) = 2^{v-2} & \text{if } p = 2, v \geq 3, \\ \varphi(p^v) = p^{v-1}(p-1) & \text{if } p > 2. \end{cases}$$

It gives the maximal order of the elements of the residues ring $\mathbb{Z}_N$. These two functions satisfy the divisibility relation $\lambda(N) \mid \varphi(N)$. There are infinitely many integers such that $\lambda(N) = \varphi(N)$, exempli gratia, $N$ = prime, twice a prime, et cetera. However, $\lambda(N) < \varphi(N)$ for almost all integers. In addition, for the integers of interest, the analysis of the equations $de - k\varphi(N) = 1$ and $de - k\lambda(N) = 1$ are almost the same.

**Theorem 25.** (i) The average order of the probability of selecting a relatively prime integer $a < N \leq X$ from the interval $[1, X)$ is

$$\sum_{1 \leq N \leq X} \frac{\varphi(N)}{N} = \frac{6}{\pi^2} X + O(\log X). \qquad (26)$$

(ii) The average order of the number of relatively prime integers $a < N \leq X$ is

$$\sum_{1 \leq N \leq X} \varphi(N) = \frac{3}{\pi^2} X^2 + O(X \log X).$$

Consult [AP, p. 71] for more details.

**Theorem 26.** (Cesaro 1881) (i) A pair of random integers $a, b \in \mathbb{N}$ are relatively prime with probability $6/\pi^2$. (ii) The average order of the number of pairs of relatively prime integers is

$$\sum_{1 \le a, b \le N} \delta(1 - \gcd(a, b)) = \frac{6}{\pi^2} N^2 + O(N \log N), \tag{27}$$

where $\delta(x) = 1$ if $x = 0$, otherwise $\delta(x) = 0$.

The proportion of relatively prime integers in short interval can be determined by mean of the Legendre sieve. The Legendre sieve is a reformulation of the sieve of Eratosthenes, circa 200 BC, in terms of number theoretical functions.

**Theorem 27.** Let $0 < X < Y \le N$ be sufficiently large integers. Then the interval $[X, Y] \subset [1, N]$ contains

$$\frac{\varphi(N)}{N}(Y - X) + O(2^{\omega(N)}), \tag{28}$$

relatively prime integers $a \in [X, Y]$, $\gcd(a, N) = 1$, where $\omega(N) = \#\{\text{ prime } p \mid N \}$.

Proof: Let $A = [X, Y] \subset [1, N]$ and let $P = \{ \text{prime } p \mid N \}$. Now apply the Legendre sieve $S(A, P) = \#\{ a \in A : \gcd(a, p) = 1, \text{and } p \in P \}$. To proceed put $Z = Y - X$. Then

$$S(A, P) = \sum_{a \in A} \sum_{d \mid a, d \mid N} \mu(d) = \sum_{d \mid N} \mu(d) \left[\frac{X}{d}\right] = X \sum_{d \mid N} \frac{\mu(d)}{d} = \frac{\varphi(N)}{N} X + O(2^{\omega(N)}). \tag{29}$$

Further, for large $N$ and $Y = N^\alpha$, $\alpha > 0$, and using the estimates $\varphi(N) = O(N/\log\log N)$ and $2^{\omega(N)} = O(2^{\log N/\log\log N})$, the proportion relatively prime integers given by

$$\frac{1}{Z}\left(\frac{\varphi(N)}{N} Z + O(2^{\omega(N)})\right) = O(1/\log\log N) \tag{30}$$

is non negligible. $\blacksquare$